\documentclass[12pt,a4paper,reqno]{amsart}

\usepackage{amsmath}
\usepackage{amsfonts}
\usepackage{amssymb}
\usepackage{amsthm}
\usepackage{amsaddr}

\usepackage[utf8]{inputenc}
\usepackage[T1]{fontenc}
\usepackage[mathscr]{eucal}
\numberwithin{equation}{section}
\usepackage[margin=2.9cm]{geometry}

\usepackage{graphicx}
\graphicspath{{./}}

\usepackage{caption}

\usepackage[hyperfootnotes=false]{hyperref}

\newcommand{\R}{\mathbb{R}}
\newcommand{\norm}[1]{\left\lVert#1\right\rVert}

\newcommand\blfootnote[1]{%
  \begingroup
  \renewcommand\thefootnote{}\footnote{#1}%
  \addtocounter{footnote}{-1}%
  \endgroup
}

\title[Mesh motion in FSI with deep operator networks]{Mesh motion in fluid-structure interaction with deep operator networks}

\author{Ottar Hellan}
\address{Simula Research Laboratory, Oslo, Norway, hellanottar@simula.no }

\begin{document}

\maketitle

\begin{abstract}\label{sec:abstract}
A mesh motion model based on deep operator networks is presented. 
The model is trained on and evaluated against a biharmonic mesh motion model on a fluid-structure interaction benchmark problem and further evaluated in a setting where biharmonic mesh motion fails.
The performance of the proposed mesh motion model is comparable to the biharmonic mesh motion on the test problems.
\end{abstract}

\section{Introduction}\label{sec:introduction}

Recently, machine learning algorithms have been developed to learn mappings between function spaces, such as the solution operators of PDEs, in a field called operator learning, with examples including the deep operator network \cite{deeponet, deeponet-comparison} and Fourier neural operator \cite{li2021fno}. The motivating applications for these models is to replace entire solvers with data-driven models, but the training of these is challenging, sometimes resulting in relative errors of several percent \cite{deeponet-comparison}.

Alternatively, operator learning can be applied to replace only \textit{certain components} of solvers, for instance targeting those steps which are based on heuristics or where classical approaches yield insufficient performance. 
One example of such components is mesh motion, where a computational reference domain is moved to fit the domain of interest. The outcome of the final solver depends on a minimal level of quality of the mesh movement.

In this work, we present a deep operator network-based mesh motion operator, trained on data from an FSI benchmark problem. 
We evaluate the learned mesh motion on this FSI benchmark, comparing it with the biharmonic mesh motion it is trained on, as well as a setting constructed such that the biharmonic mesh motion fails.

\section{Mesh motion}\label{sec:mesh_motion}

The problem of mesh motion is prominent in a number of applications involving a change-of-coordinates between reference and target domains. 
In these cases, a mesh motion operator is used to move a mesh of the reference domain to triangulate the target domain. 
Examples include fluid-structure interaction (FSI) in the arbitrary Lagrangian-Eulerian (ALE) formulation \cite{wick2011, hron2006monolithic}, where the mesh of the fluid domain is moved to fit a time-changing solid domain, or in shape optimization via the method of mappings, where a reference mesh is moved to fit the desired shape \cite{haubner2021shapeopt}. 
The change-in-coordinates given by the mesh motion are typically based on PDEs, with Dirichlet boundary conditions representing the geometry of the target domain \cite{shamanskiy2020}, but recently data-driven approaches have been suggested \cite{aygun2023mmpinns, song2022m2n, learnext}.

Let the domain $\Omega$ be a bounded subset of $\R^d$, $d \in \{2,3\}$. 
We define a mesh motion operator as a mapping $g \mapsto u$, 
where $g: \partial \Omega \to \R^d$ represents the deformation of the boundary $\partial\Omega$ and  
the resulting deformation field $u: \Omega \to \R^d$ defines the change-of-coordinates
\begin{equation}\label{eq:mm_changecoords}
    \chi: \Omega \to \R^d, \quad x \mapsto x + u(x),
\end{equation}
with the boundary condition $u|_{\partial\Omega} = g$. 
For this change-of-coordinates to be well defined and suitable for our computations, it is necessary that it is bijective, that $u|_{\partial\Omega} = g$, and that 
\begin{equation}\label{eq:det_grad_chi}
    J = \mathrm{det}(\nabla \chi) > 0 \text{ in } \Omega.
\end{equation}
In addition to these constraints the transformation of $\Omega$ should produce a suitable domain for the computational problem, for instance having a Lipschitz boundary, and the resulting mesh should be of high quality. Nevertheless, there is considerable freedom in the choice of $u$, allowing the definition of many different mesh motion operators.

The simplest mesh motion operator is the harmonic mesh motion \cite{wick2011}, which is defined by the solution operator of the Laplace equation
\begin{equation} \label{eq:harmonic_mm}
-\Delta u = 0 \text{ in } \Omega, \quad u = g \text{ on } \partial\Omega, 
\end{equation}
where $\Delta \cdot$ is the vector-valued Laplacian $(\Delta u)_i = \sum_{j=1}^d \partial^2 u_i / \partial x_j^2$ for $j = 1, ..., d$.
Another common mesh motion operator, which can deal with larger deformations than harmonic mesh motion, is the biharmonic mesh motion \cite{wick2011}, one variant of which is stated as 
\begin{equation}\label{eq:biharmonic_mm}
    \Delta^2 u = 0 \text{ in } \Omega, \quad u = g \text{ on } \partial\Omega, \quad \nabla u \cdot \mathbf{n} = 0 \text{ on } \partial\Omega.
\end{equation}
Other versions of biharmonic mesh motion exist which use different sets of boundary conditions \cite{wick2011}.

\section{Deep Operator Networks}\label{sec:deeponet}

The deep operator network (DeepONet) is a neural network model that is used to learn operators, that is, mappings between function spaces \cite{deeponet}. We wish to learn mesh motion operators, and thus we map boundary deformations $g \in \mathcal{G}$ to the deformation map $u$ that defines the change-of-coordinates \eqref{eq:mm_changecoords}.

For DeepONets with scalar output functions, we evaluate the output function $\mathscr{D}(g)$ of the learned operator at a point $x$ by the inner product of two separate networks,  
\begin{equation}\label{eq:deeponet_def}
    \mathscr{D}(g)(x) = \mathscr{B}(\mathscr{E}_\mathscr{B}(g)) \cdot \mathscr{T}(x).
\end{equation}
The neural networks $\mathscr{B}: \R^{b_\mathrm{in}} \to \R^p$ and $\mathscr{T}: \R^{t_\mathrm{in}} \to \R^p$ are known as the \textit{branch} and \textit{trunk networks} respectively and can be of any given architecture, for instance multi-layer perceptrons (MLPs), convolutional neural networks, or graph neural networks \cite{deeponet-comparison}. 
The mapping $\mathscr{E}_{\mathscr{B}}: \mathcal{G} \to \R^{b_\mathrm{in}} $ encodes the input function $g$ to discrete, finite-dimensional values that neural networks can process. Typically, $\mathscr{E}_\mathscr{B}$ is the evaluation of $g$ at a fixed set of locations called \textit{sensors}, but can also be projection onto a predefined basis \cite{deeponet-comparison}. 
To work with vector-valued output functions, we use the second of the four approaches presented in \cite{deeponet-comparison}. The first $p / d$ components of $\mathscr{B}(\mathscr{E}_\mathscr{B}(g))$ and $\mathscr{T}(x)$ are used to construct $u(x)_1$, the first component of $\mathscr{D}(g)(x)$, the next $p / d$ components for $u(x)_2$, and so on, producing an output vector in $\R^d$.

In our DeepONet models, we let $\mathscr{E}_\mathscr{B}$ be the concatenation of $g$ evaluated at a number of sensors placed appropriately in the domain of $g$, and let $\mathscr{B}$ and $\mathscr{T}$ be MLPs. The output of an MLP is defined by 
\begin{subequations}\label{eq:mlp_def}
\begin{align}
    & \mathcal{N}(z) = W_L z_{L-1} + b_L, \\
    & z_l = \sigma(W_l z_{l-1} + b_l), \text{ for } l = 2, ..., L-1, \quad z_1 = z,
\end{align}
\end{subequations}
with $\sigma$ a chosen activation function.  
The matrices $W_2 \in \R^{w\times\mathrm{dim}(z)}$, $W_l \in \R^{w \times w}$, $W_L \in \R^{\mathrm{dim}(\mathcal{N}(z)) \times w}$ and vectors $b_l$ are trainable parameters of the network.
The \textit{depth} of the network is $L$ and the \textit{width} of the network is $w$.

For mesh motion, it is important to satisfy the Dirichlet boundary constraints exactly. 
If this is not the case, the geometry will be represented incorrectly, causing numerical artifacts. We therefore utilize a hard constraint boundary condition technique to ensure the DeepONets produce a valid mesh motion operator. Following \cite{mcfall2009hardconstraints, deeponet-comparison}, 
we define the evaluation at $x$ of our DeepONet mesh motion model $\mathscr{U}$ as 
\begin{equation}\label{eq:corrected_deeponet}
   \mathscr{U}(g)(x) = h(g)(x) + l(x) \, \mathcal{D}(g)(x),
\end{equation}
where $h(g)$ is a function satisfying the Dirichlet boundary conditions exactly and $l: \overline\Omega \to \R$ is a function satisfying $l = 0$ on $\partial\Omega$ and $0 < l(x) < 1$ for $x \in \Omega$.

\section{FSI benchmark problem}\label{sec:fsi_problem}

\begin{figure}
    \centering
    \includegraphics[width=0.75\linewidth]{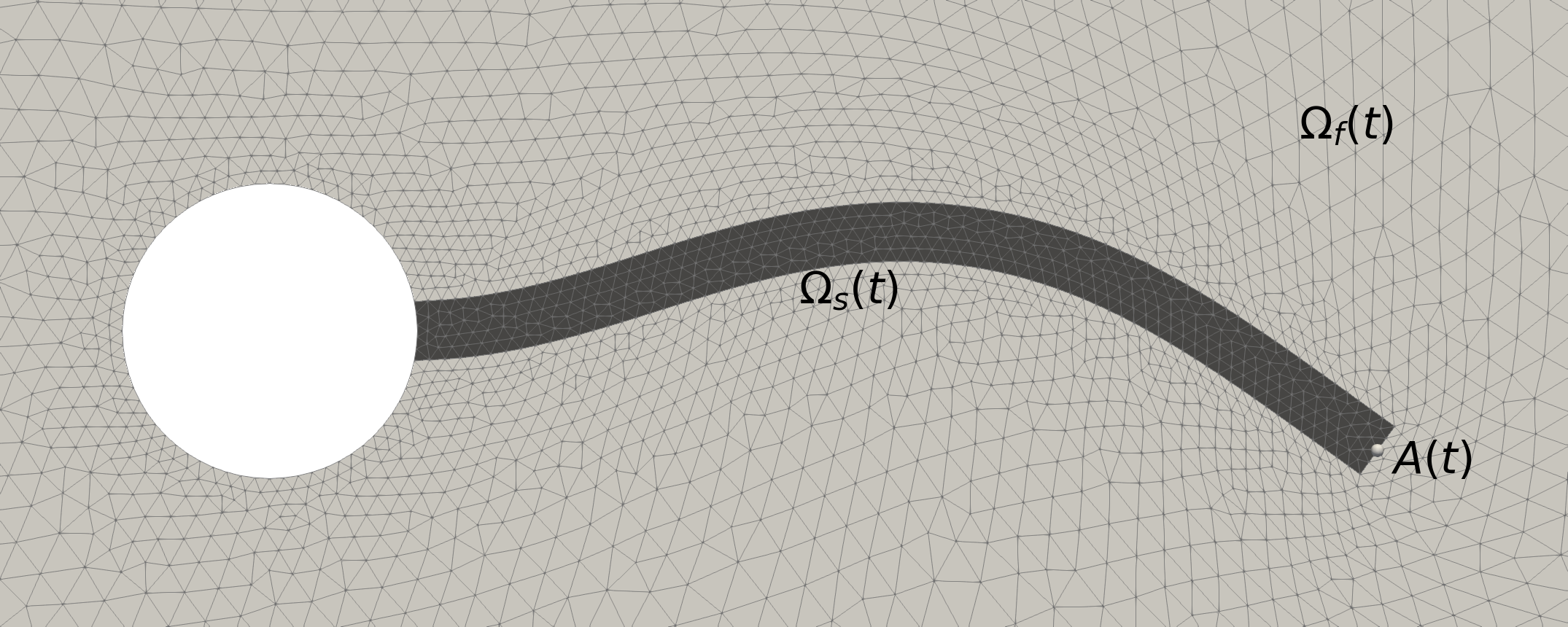}
    \caption{Partial view of the FSI2 benchmark problem's geometry, centered on the submerged solid. The displacement of the particle $A(t)$ is a reported quantity in the benchmark. }
    \label{fig:fsi2_domain_fig}
\end{figure}

To train and evaluate our DeepONet mesh motion model, we implement the FSI2 benchmark problem \cite{Turek_fsi_benchmark}, simulating the interaction of an incompressible Navier-Stokes fluid with an obstacle comprised of a rigid cylinder with an elastic solid made of compressible St. Venant-Kirchhoff material attached behind it, see Figure \ref{fig:fsi2_domain_fig}. The problem setup results in periodic oscillations of the elastic solid.

The governing equations of the solid and fluid, in the ALE formulation on a fixed reference domain $\Omega = \Omega_s \cup \Omega_f$, are given as 
\begin{subequations}\label{eq:fsi_monolithic}
\begin{align}
    \rho_s \partial^2 u_s / \partial t^2 = \mathrm{div} (J \sigma_s F^{-T}) \text{ in } \Omega_s, \\
    \rho_f \partial v_f / \partial t + \rho_f (\nabla v_f) F^{-1} (v_f - \partial u / \partial t ) = J^{-1} \mathrm{div} (J \sigma_f F^{-T})  \text{ in } \Omega_f,  \\
    \mathrm{div} (J v_f F^{-T}) = 0 \text{ in } \Omega_f, \\
    \sigma_s = J^{-1} F \left( \lambda_s (\mathrm{tr}\, E) \mathrm{I} + 2\mu_s E \right) F^T, \quad \sigma_f = -p_f \mathrm{I} + \rho_f \nu_f (\nabla v_f + \nabla v_f^T),
\end{align}
\end{subequations}
where, $F = \mathrm{I} + \nabla u$, $J = \mathrm{det}(F)$, and $E = (F^T F - \mathrm{I})/2$ \cite{Turek_fsi_benchmark, hron2006monolithic}. 
Here, $\rho$ is the density of the material, $\nu_f$ is the viscosity of the fluid, and $\lambda_s$, $\mu_s$ are the Lamé parameters of the solid.
To complete the formulation, there are appropriate initial, boundary, and coupling conditions, in particular enforcing the stress and mass balance on the fluid-solid interface. 
In $\Omega_s$, $u$ represents the deformation of the solid, and in $\Omega_f$, $u$ represents the deformation field that defines the ALE change-of-coordinates \eqref{eq:mm_changecoords}. 
The ALE formulation is introduced to reconcile the natural formulation of the solid problem in a Lagrangian frame and the fluid problem in an Eulerian frame \cite{donea2004ale, hron2006monolithic}.

A detailed description of a solver for \eqref{eq:fsi_monolithic} of monolithic type is given in \cite{hron2006monolithic} and the performance of different mesh motion operators in monolithic ALE-FSI solvers is investigated in \cite{wick2011}. 
We solve \eqref{eq:fsi_monolithic} using a monolithic solver from a previous work \cite{learnext}, which solves the resultant nonlinear system in several steps, solving for $u$ in $\Omega_f$ separately and the solid and fluid equations together, allowing us to easily change the mesh motion model used in our tests. 
The solver uses second order Lagrange finite elements for $u$ and $v$ and linear Lagrange elements for $p$ over a triangular mesh. 
We use the FSI2 benchmark problem and this solver to generate our dataset and to evaluate the trained DeepONet mesh motion.

\section{Network details}\label{sec:network}

\subsection{Training process}\label{sec:training_process}

We train our DeepONet models in a supervised fashion, with biharmonic mesh motion as target data. 
Using our FSI solver with biharmonic mesh motion, we simulate one period of the oscillatory motion in the FSI2 benchmark, after a 15 second period to allow the oscillations to fully develop, producing 207 different snapshots $u_\mathrm{bih}^k$ that make up our dataset. We randomly split this dataset into 70\% training set and 30\% validation set.

The cost function used is
\begin{equation}
    \mathcal{J} = \frac{1}{K} \sum_{k} \sum_{i} \frac{ \left( h(g^k)(x_i) + l(x_i) \, \mathcal{D}(g^k)(x_i) - u_\mathrm{bih}^k(x_i) \right)^2 }{ \left( h(g^k)(x_i) - u_\mathrm{bih}^k(x_i) \right)^2 + \varepsilon}.
\end{equation}
Here, $K$ is the number of snapshots in the batch, $g^k$ is the boundary deformation of snapshot $k$ and $u_\mathrm{bih}^k$ is the solution of \eqref{eq:biharmonic_mm} with boundary condition $g^k$.
The evaluation points $x_i$ are the mesh vertices of our mesh of the FSI2 geometry. 
The parameter $\varepsilon > 0$ adds stability, since $h(g^k)(x_i) - u^k_\mathrm{bih}(x_i)$ is close to zero for some $i$ and $k$.

Our models and training are implemented in the \texttt{PyTorch} \cite{paszke2019pytorch} deep learning-framework. We use the built-in Adam optimizer \cite{kingma2017adam}, 
with learning rate $10^{-5}$ and otherwise default parameters, 
and the built-in \texttt{ReduceLROnPlateau} learning rate scheduler, 
with the reduction factor set to $0.5$ and otherwise default parameters, for $40000$ epochs.

\subsection{Hyperparameters}\label{sec:hyperparameters}

Following \cite{learnext}, we define our DeepONet mesh motion model by \eqref{eq:corrected_deeponet}, with $h(g)$ chosen as the harmonic extension \eqref{eq:harmonic_mm} of the boundary deformation $g$. Harmonic mesh motion is suitable for small boundary deformations and satisfies the boundary deformation exactly. 
It is also a well studied PDE, with fast, scalable solvers \cite{gholami2016poisson}.

The function $l: \overline\Omega \to \R$ in \eqref{eq:corrected_deeponet} is determined by solving a Poisson problem
\begin{equation}
    -\Delta \tilde l(x) = f(x) \text{ in } \Omega, \quad \tilde l(x) = 0 \text{ on } \partial\Omega,
\end{equation}
with $f(x) = 2(x+1)(1-x)\mathrm{exp}(-3.5 x^7) + 0.1$ specifically chosen to weight areas of the FSI2 benchmark problem's domain where harmonic mesh motion is inadequate \cite{learnext}. 
From classical PDE theory, it holds that $0 < \tilde l$ in $\Omega$ and $\tilde l \in C^2(\Omega) \cap C(\overline\Omega)$. We then take $l = \tilde{l} / \norm{\tilde l}_\mathrm{sup}$, such that $0 < l < 1$ in $\Omega$. With $f \equiv 1$ the same properties hold for our and also for many other relevant geometries \cite{learnext}.

We choose the sensor locations of $\mathscr{E}_\mathscr{B}$ encoding $g$ to be the mesh vertices of the solid domain's boundary, as this encodes the geometry of the deformed solid in our simulations. On our mesh, this entails evaluating $g$ at $206$ vertices, resulting in $b_\mathrm{in} = 412$-dimensional input-vectors to the branch network. The trunk network takes inputs of size $t_\mathrm{in} = d = 2$.

With our DeepONet mesh motion model and training process defined, we search for appropriate values of the model's hyperparameters, to obtain accurate results on our test problems. 
We choose the hyperparameters of our DeepONet by performing a grid search over a selection of hyperparameter values, leaving the rest fixed at values we found were robust choices for many different architectures in the exploratory phase of the research. 
The branch and trunk network architectures are constrained to be equal, except for the size of the input layer, and we vary the depth $d$, width $w$, and output size $p$ of the two MLPs.

To evaluate our mesh motion models, we use the scaled Jacobian mesh quality measure on the deformed meshes, computed using the visualization toolkit \cite{vtk2006book}. This mesh quality measure is scale invariant and gives the value of 1 for equilateral triangles and goes to zero as a cell degenerates.
Since the networks are trained in a supervised fashion against the biharmonic extension, we assume their performance will be worse than the biharmonic extension and select the DeepONet with the lowest drop in quality. 
We compute for each snapshot $k$ the difference in minimal cell value of mesh quality between the biharmonic and DeepONet mesh motion, and select the DeepONet with the lowest such difference.

The hyperparameter values used in the grid search are listed in Table \ref{tab:hyperparameter_study_values}. All hyperparameter combinations are tested twice, using the random seeds 0 and 1, to account for random initialization. All other network and training parameters are unchanged between runs. The number of network parameters ranges from $160\,576$ to $3\,430\,528$. 

\begin{table}
    \centering
    \begin{tabular}{c | c c c c c}
        \hline
        hyperparameter & depth ($d$) & \quad & width ($w$) & \quad & output size ($p$) \\
        values & $\left\{4, 5, 6, 7\right\}$ & \quad & $\left\{128, 256, 512\right\}$ & \quad & $\left\{32, 64\right\}$\\
        \hline
    \end{tabular}
    \caption{Hyperparameter values used grid search hyperparameter optimization.}
    \label{tab:hyperparameter_study_values}
\end{table}

\begin{figure}
    \centering
    \includegraphics[width=0.49\textwidth]{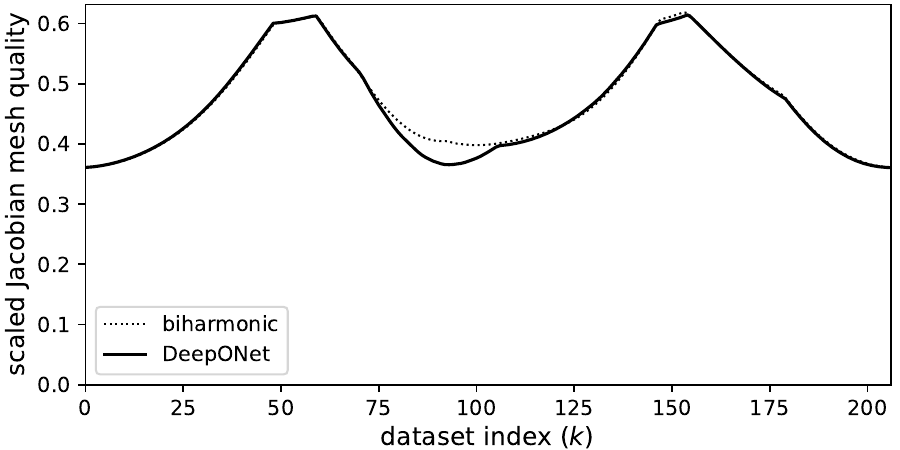}
    \caption{Minimum cell value of scaled Jacobian mesh quality measure over the dataset with biharmonic extension and best trained DeepONet ($d=7$, $w=512$, $p=32$, seed 1).}
    \label{fig:mesh_quality_hyperparameters}
\end{figure}

Figure \ref{fig:mesh_quality_hyperparameters} shows the minimal cell value of the scaled Jacobian mesh quality for the biharmonic extension and the best trained DeepONet extension over the dataset. 
The resulting extension is of high quality and the worst meshes from the DeepONet mesh motion are not poorer than the worst meshes from the biharmonic mesh motion in this dataset. 
The chosen architecture is one of the largest, with both depth $d$ and width $w$ being the largest possible values, but the networks' output size $p$ is the lower of the two possible values. 
We are testing on data included in the training set, but with a test metric different from the cost function.

To investigate the robustness of the model and training process, we investigate the dependence of random initialization on the quality of the mesh motion. We train 20 different random initializations of the DeepONet model selected in the hyperparameter search. All parameters are kept constant, except for the random seeds, which are chosen from $\{0, ..., 19\}$.

\begin{figure}
    \includegraphics[width=0.49\textwidth]{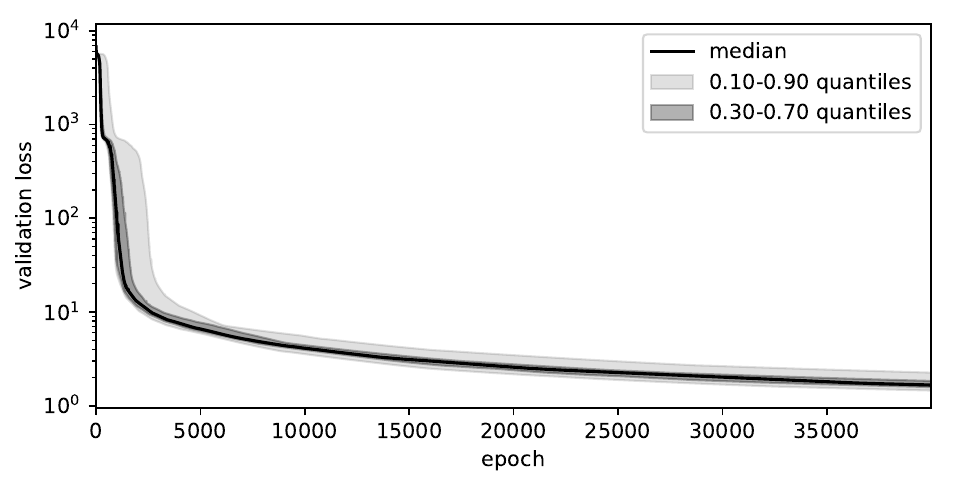}	\includegraphics[width=0.49\textwidth]{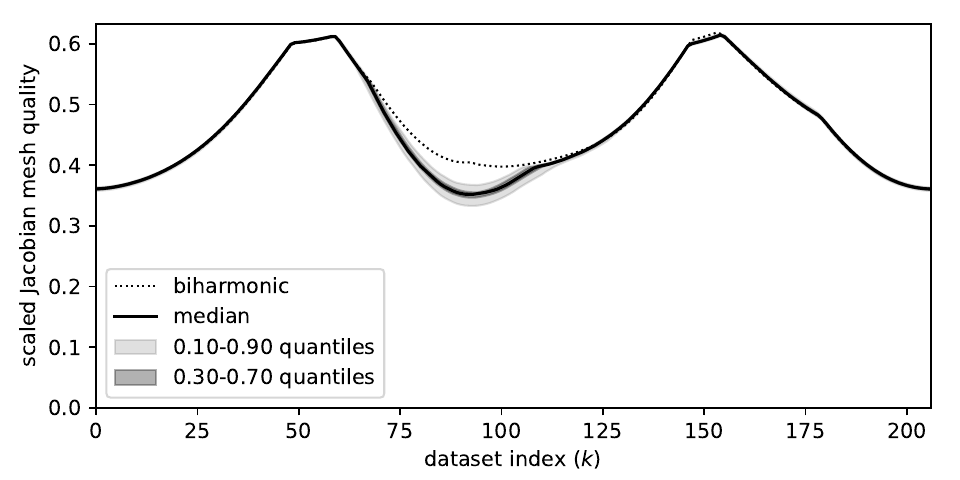}
    \caption{Quantiles of validation loss history DeepONet training (left) and resulting scaled Jacobian mesh quality over FSI2 dataset (right) for 20 random initializations of the best performing hyperparameters found in grid search.}
    \label{fig:random_initialization_results}
\end{figure}

In Figure \ref{fig:random_initialization_results}, we report quantiles of the validation loss history and the minimum cell value of scaled Jacobian mesh quality at dataset index $k$ for the fully trained DeepONet mesh motion models. The results indicate that the model and training process result in a high quality mesh motion model, with $0.1$-$0.9$ quantiles of minimal cell value scaled Jacobian mesh quality almost on par with biharmonic mesh motion. 
Of the 20 random initializations, 19 of them give results similar to those shown in Figure \ref{fig:random_initialization_results}, while one initialization stagnated in a local minimum almost immediately, after approximately 10 epochs.

\section{FSI test}\label{sec:fsi_test}

To evaluate the DeepONet mesh motion, we use it as the mesh motion operator in the FSI2 benchmark problem with our monolithic solver. 
We evaluate the DeepONet mesh motion on a GPU, while the coupled fluid and solid equations are solved in serial on the CPU.
The mesh motion is performed by computing $\mathcal{D}(g)(x_i)$ for every degree of freedom-location $x_i$ of the finite element space of $u$, not only the mesh vertices, thereafter inserting the values in the basis coefficient array of the finite element function.
The mesh used in the FSI simulations is the same on the solid boundary as in the training process, up to ordering, so to encode the input function $g$, we access the finite element basis array using a pre-built vector of indices. If the mesh vertices and sensors do not align, this could be done similarly with a pre-built, sparse interpolation matrix. The functions $h(g)$ and $l$ are computed using finite elements, with the same function space as $u$ is represented in.

We report some quantities of interest from the FSI2 benchmark to validate the DeepONet mesh motion. 
The simulation is run for 15 seconds, allowing the periodic oscillations of the solid to fully develop, using the DeepONet mesh motion. 
Thereafter, we compute for the period $[15\,\mathrm{s},16\,\mathrm{s}]$ the drag and lift forces applied to the solid by the fluid flow and the minimal value in $\Omega$ of $J$, see \eqref{eq:det_grad_chi}. 
We also report the $y$-displacement $A_2$ of a point in the solid, see Figure \ref{fig:fsi2_domain_fig}. 
The same values are also computed using biharmonic mesh motion, initialized at $t=15\,\mathrm{s}$ with the same state as the DeepONet mesh motion.

\begin{figure}
    \includegraphics[width=0.49\textwidth]{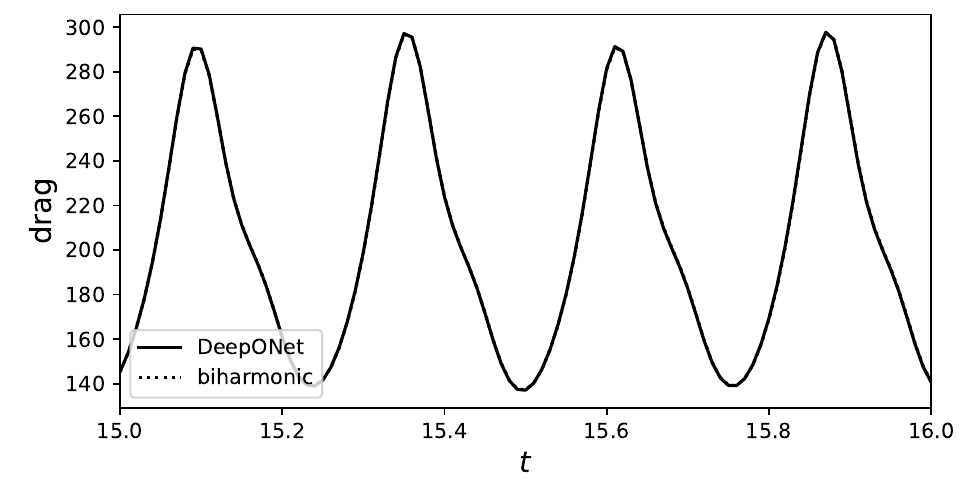}	\includegraphics[width=0.49\textwidth]{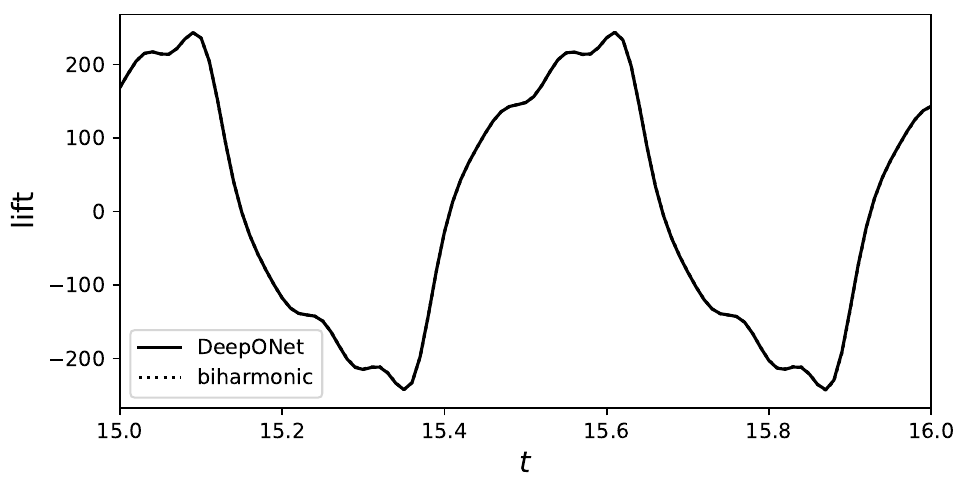} \\
    \includegraphics[width=0.49\textwidth]{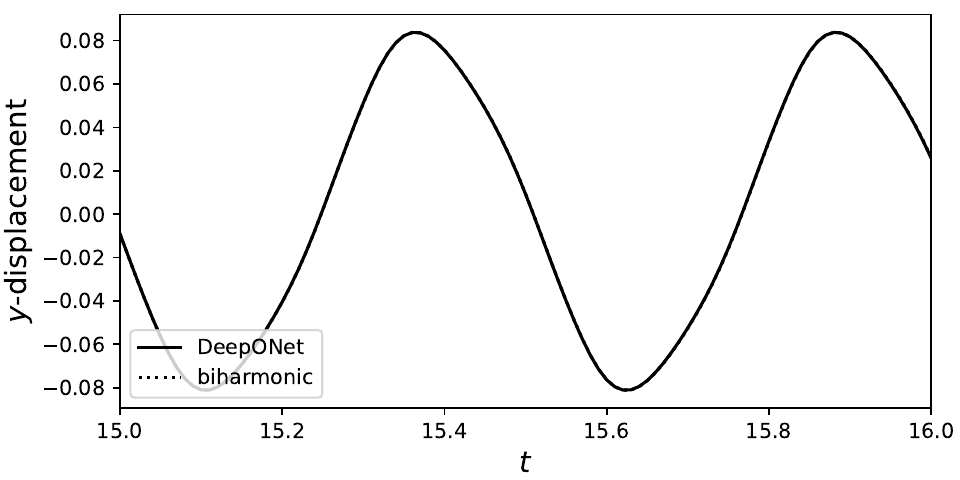}	
    \includegraphics[width=0.49\textwidth]{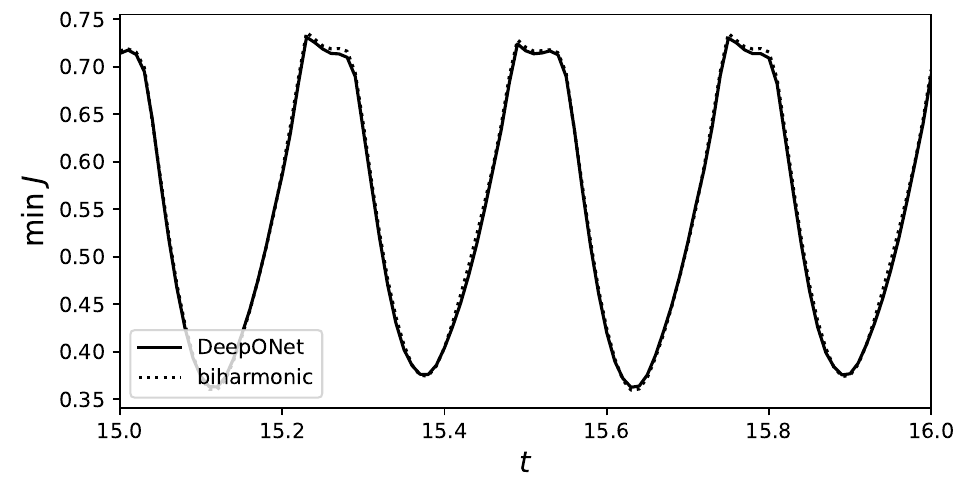}	
    \caption{Drag (top left), lift (top right), $y$-displacement $A_2$ of point on tip of solid (bottom left), and minimum value in $\Omega$ of $J$ (bottom right) in FSI2 benchmark problem with biharmonic and DeepONet mesh motion.}
    \label{fig:FSI_results}
\end{figure}

Figure \ref{fig:FSI_results} shows the computed quantities of interest. 
The drag, lift, and $y$-dis-placement values are identical for the two mesh motion models and match the reference values of \cite{Turek_fsi_benchmark} well. This indicates the DeepONet mesh motion is suitable for FSI computations in ALE formulation. The value of $\min J$ is slightly lower at certain times for the DeepONet mesh motion, indicating a slightly lower quality mesh motion than the biharmonic. This is in line with the results in Figure \ref{fig:random_initialization_results}.

\section{Gravity driven deformation test}\label{sec:grav_test}

To further explore the performance of the DeepONet mesh motion model, we apply it to a test problem \cite[sec. 4.1]{shamanskiy2020} constructed to push biharmonic mesh motion to degeneracy. With the same geometry, material model, and parameters as in the FSI2 benchmark, the solid is subjected to a uniform gravitational load $(0, f_g)$ and the solid deformation at its maximum deformation is extended to the fluid domain with different mesh motion models. The solid deformations are computed using implicit Euler time discretization and second order Lagrange finite elements in space.

\begin{figure}
    \centering
    \includegraphics[width=0.9\linewidth]{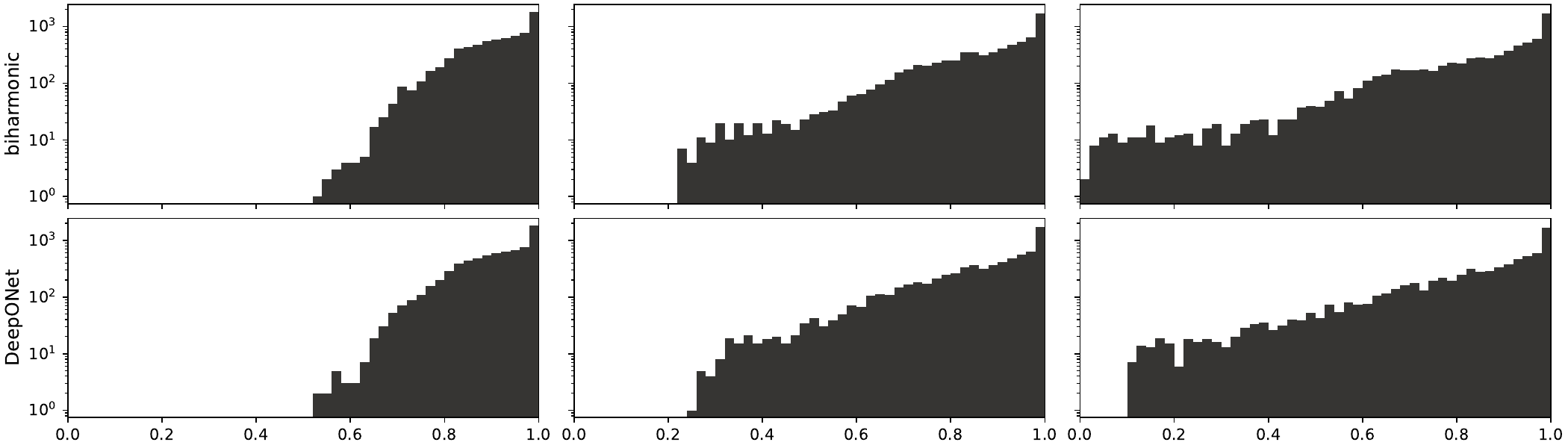}
    \caption{Histograms of scaled Jacobian mesh quality for biharmonic and DeepONet mesh motion models with solid deformations caused by uniform gravitational load $(0, f_g)$, $f_g \in \{ 1, 2, 2.5 \}$. }
    \label{fig:grav_hist}
\end{figure}

Figure \ref{fig:grav_hist} shows histograms of the scaled Jacobian mesh quality produced by the biharmonic and DeepONet mesh motions with the maximal solid deformations caused by the gravitational loads $f_g \in \{1, 2, 2.5\}$ as boundary conditions. 
The biharmonic mesh motion handles the first two deformations well, but breaks down for the third due to collapsing cells, with mesh quality going to zero. Although the DeepONet mesh motion is trained to learn the biharmonic mesh motion, it is able to handle even the largest deformation.

\section{Conclusion and future work}\label{sec:conclusion}

We have presented a DeepONet-based mesh motion model trained on biharmonic mesh motion data, with performance on the test problems comparable to the biharmonic mesh motion it is trained on. 
The training data comes from the same problem as one of the test problems and future work should therefore investigate training the model without access to this data. 
Training the network in an unsupervised fashion and applying the method to problem settings where current mesh motion operators struggle would be especially interesting. 
Other work on the efficiency and scalability of the method can also be considered.

\bibliographystyle{abbrv}
\bibliography{references}

\blfootnote{Ottar Hellan acknowledges support from the Research Council of Norway, grant 303362.}

\end{document}